# Intervalles de confiance pour une proportion : lesquels doit-on enseigner ?

*Jean-Christophe Turlot, IUT STID /UPPA/Paujean-christophe.@univ-pau.fr*
*Jean-François Petiot, IUT STID /UBS/Vannesjean-francois.petiot@univ-ubs.fr*

**Résumé.** Les intervalles de confiance pour une proportion les plus usuellement enseignés sont l'intervalle de Wald (Ws) et l'intervalle de Cloper-Pearson (CP), pour la simplicité de leur définition. Or, la probabilité réelle de recouvrement du paramètre $p$ est erratique, souvent bien éloignée de la probabilité nominale souhaitée de $1-\alpha$. Si d'autres intervalles sont bien préférables à ceux-ci, leur expression est généralement complexe et peu interprétable. Or, par une modification simple des définitions des intervalles Ws et AC, on en obtient d'autres avec des probabilités de recouvrement bien plus satisfaisantes. Ce sont les intervalles d'Agresti-Coull et Mid-P que nous présentons ici. Nous les recommandons vivement dans le cadre d'enseignements de base en statistique.

**Mots clés.** Intervalle de confiance pour une proportion, probabilité réelle de recouvrement, développement d'Edgeworth, précision d'un intervalle.

**Abstract.** The most frequently taught confidence intervals for a proportion are the classical Wald (Ws) and the Clopper-Pearson (CP) ones because of the simplicity of their definition. However, their actual coverage probability of the parameter p is erratic, often quite far from the nominal probability $1-\alpha$ which is aimed at. Other confidence intervals are clearly preferable to the former, but their expression is generally complex and they are difficult to interpret. But nevertheless, through a simple modification of the definition of the Ws and CP intervals, we obtain some confidence intervals with much better coverage probabilities. Namely, these confidence intervals are the Agresti-Coull and Mid-P intervals that we present here. We highly recommend them in a basic Statistics course.

**Keywords:** Confidence interval for a proportion, actual coverage probability, Edgeworth expansion, accuracy of an interval.

## Introduction

La notion d'intervalle de confiance est très souvent mise en application dans tout champ d'activité faisant appel à la statistique inférentielle. Dans le cas de l'estimation d'une proportion $p$, la plupart des ouvrages généraux de statistique présentent l'intervalle standard de Wald ($Ws$) qui présente l'avantage d'avoir une expression simple : $IC_{Ws} = \left[\hat{p} \pm \kappa\sqrt{\hat{p}(1-\hat{p})/n}\right]$, où la valeur de $\kappa$ est choisie selon le niveau de confiance souhaité (typiquement 95%) en faisant appel à l'approximation normale de $\hat{p}$ (l'estimateur naturel de $p$) ce qui suppose que $n$ n'est pas trop petit. Un autre intervalle classiquement utilisé en contrôle de qualité dans l'industrie est obtenu par l'inversion directe de l'intervalle de fluctuation de $n\hat{p}$ qui suit une loi $B(n,p)$ appelé intervalle « exact » de Clopper-Pearson $(CP)$.

On montre qu'aucun de ces deux intervalles n'est satisfaisant au sens où la probabilité réelle de recouvrement de $p$ sont généralement bien loin de leur valeur nominale de $1-\alpha$ lorsque la taille de l'échantillon est modérée $(n < 100)$. Si la présence d'oscillations dans la probabilité vraie de recouvrement est inhérente à la nature qualitative de la loi binomiale (ce qui empêche d'obtenir des intervalles de niveau exact), ces deux intervalles présentent en plus un biais de recouvrement trop important: $IC_{WS}$ est trop laxiste en ce sens que la vraie probabilité de recouvrement de $p$ est inférieure à la valeur nominale (hors oscillations), alors que $IC_{CP}$ est trop conservatif car la probabilité de recouvrement de $p$ dépasse trop largement la spécification, ce qui rend cet intervalle trop imprécis. Ces biais peuvent être réduits par des modifications simples de ces deux intervalles qui sont les plus couramment présentés dans les enseignements.

## 1. L'intervalle d'Agresti-Coull (AC)

L'intervalle de confiance d'Agresti-Coull (1998) est un compromis entre l'intervalle de Wald standard et l'intervalle de Wilson. Il adopte la simplicité d'expression du premier et approche la qualité de recouvrement du second. L'intervalle de Wald standard est fondé sur l'approximation $W_n = \dfrac{\sqrt{n}(\hat{p} - p)}{\sqrt{\hat{p}\hat{q}}} \xrightarrow{L} N(0,1)$. En utilisant cette loi limite, on obtient l'intervalle de Wald $IC_{WS}$ :

$IC_{WS} = \hat{p} \pm \kappa\sqrt{\hat{p}(1-\hat{p})/n}$ . Cet intervalle n'est pas satisfaisant, car la probabilité de recouvrement du paramètre $p$ est erratique ; en particulier, outre de fortes oscillations liées au caractère discret de la loi binomiale lorsque $p < 0.2$ ou $p > 0.8$, il présente une probabilité réelle de recouvrement généralement bien inférieure à sa valeur nominale : $IC_{WS}$ est trop laxiste [Fig.1]. Pourquoi ?

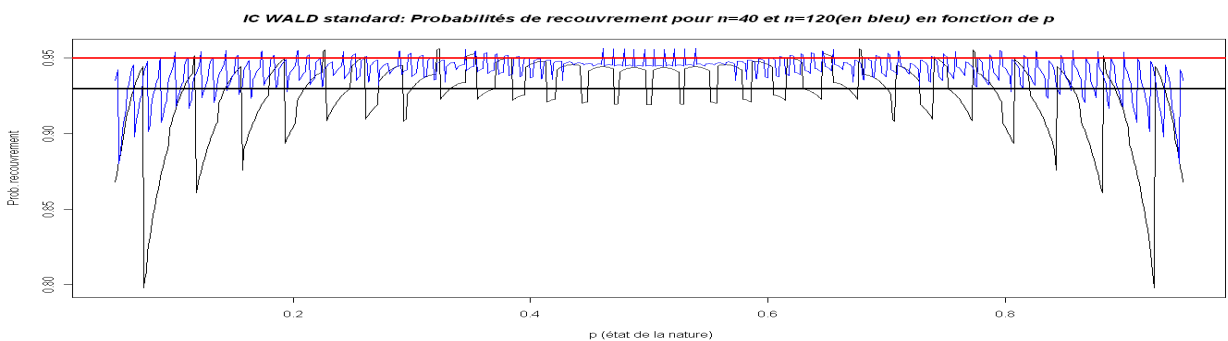

**Fig.1 : Intervalle de Wald standard : probabilités exactes de recouvrement de $p$ au niveau nominal de $1 - \alpha = 0.95$ ($n = 40$ en noir, $n = 120$ en bleu).**

Le centrage de cet intervalle en $\hat{p}$ n'est pas une bonne idée. En effet, pour $n$ modéré, l'espérance de $W_n$ est relativement éloignée de la valeur 0. Par un simple développement limité, on montre que $E(W_n) = \dfrac{p - 1/2}{\sqrt{npq}} + o(1/n)$. Il en résulte que la probabilité réelle de l'intervalle de fluctuation de $p$ est inférieure à sa valeur nominale, cela se reporte sur l'intervalle de confiance obtenu par inversion. L'idée est donc de recentrer l'intervalle fondé sur $W_n$ par son espérance approchée à l'ordre $o(1/n)$. L'intervalle de Wilson $IC_W$ est fondé sur l'approximation normale de $Z_n$ :

$Z_n = \sqrt{n(\hat{p}-p)}\big/\sqrt{p(1-p)} \approx N(0,1)$ répond à un centrage correct puisque $E(Z_n)=0, \forall n$. Il correspond à l'inversion de la famille de régions d'acceptation des tests bilatéraux $H_0 : p = p_0$, $p_0 \in (0,1)$ fondés sur la statistique $\hat{p}$. Cependant, l'expression de cet intervalle de confiance obtenu par inversion est complexe et difficilement directement interprétable.

La proposition d'Agresti-Coull : ces deux auteurs ont observé que le centre de cet intervalle est une moyenne pondérée de l'estimateur $\hat{p}$ et $1/2$ : $\tilde{p}_W = \hat{p}\big(n/(n+\kappa^2)\big) + 0.5\big(\kappa^2/(n+\kappa^2)\big)$. Cet estimateur tend à recentrer l'estimateur naturel vers 0.5 (rétrécissement). Comme la valeur $\kappa$ est très proche de 2 si l'on choisit une probabilité de recouvrement $(1-\alpha)=0.95$, Agresti et Coull ont proposé de rajouter artificiellement à l'échantillon 4 observations répondant à un tirage déterministe : la règle +2 succès, +2 échecs. On a alors : $\tilde{p}_{AC} = \hat{p}\big(n/(n+4)\big) + 0.5\big(4/(n+4)\big) \cong \tilde{p}_W$

De plus, si l'on complète l'échantillon $(X_1, X_2,..., X_n)$ par $(X_{n+1},...,X_{n+4}) \to B(1/2)$, on peut vérifier par un calcul direct que deux fois l'écart type de cet estimateur correspond à la demi-longueur de l'IC de Wilson. En remplaçant les quatre variables aléatoires rajoutées à l'échantillon par deux succès et deux échecs (de manière déterministe), l'intervalle d'Agresti Coull [Fig.2] a une expression très simple : $IC_{Ws} = \tilde{p}_{AC} \pm \kappa\sqrt{\tilde{p}_{AC}(1-\tilde{p}_{AC})/\tilde{n}}$ avec $\tilde{n} = n+4$.

Il est donc facile à présenter et très proche de l'intervalle de Wilson qui a un très faible biais de recouvrement (hors oscillations) en même temps qu'il est précis, comme on le montre maintenant.

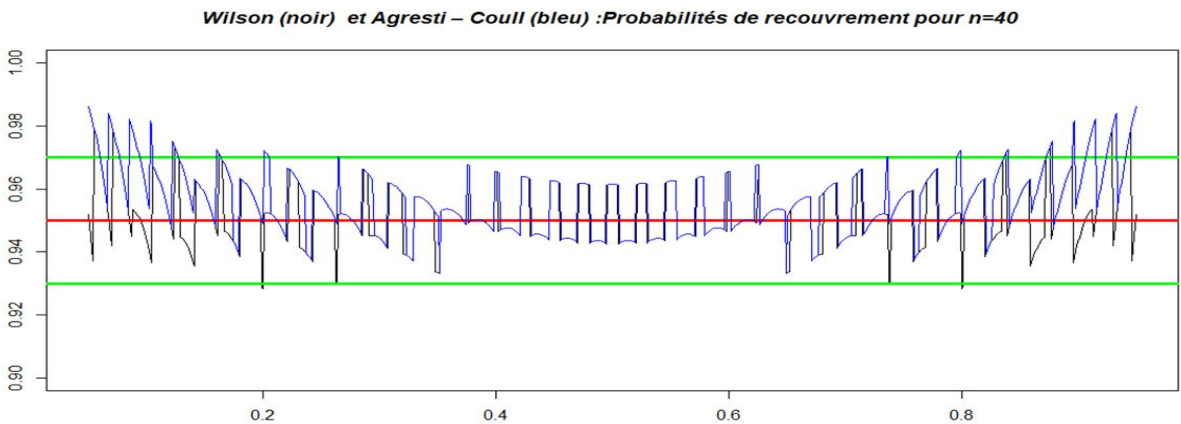

**Fig.2 : Intervalle d'Agresti-Coull : probabilités exactes de recouvrement de *p* au niveau nominal de** $1-\alpha = 0.95$ **(** $n = 40$ **en noir,** $n = 120$ **en bleu).**

L'approximation d'Edgeworth de la probabilité de recouvrement de *p* fournit une approximation analytique très précise de cette probabilité en fonction de $(p,n)$. Soit $IC_*(p,n) = [l_*(\hat{p},n), u_*(\hat{p},n)]$ un intervalle de confiance quelconque, et $C_*(p,n) = Pr_p(p \in [l_*(\hat{p},n), u_*(\hat{p},n)])$ la probabilité de recouvrement de *p* associée. Un développement d'Edgeworth à l'ordre 2 de $Pr(Z_n \le z)$ où $z = z(n)$ [Brown et al. 2002] fournit une approximation analytique de $C_*(p,n)$ de la forme suivante :

$$C_*(p,n) = (1-\alpha) + \text{oscil. en } O(n^{-1/2}) + \text{biais en } O(n^{-1}) + \text{oscil. en } O(n^{-1}) + O(n^{-3/2})$$

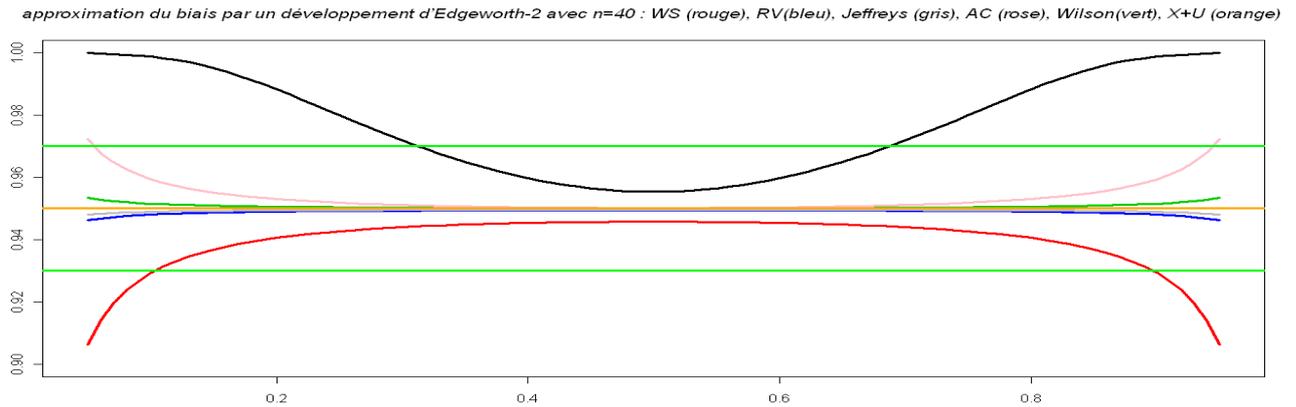

Fig.3 : Approximation du biais par un développement d'Edgeworth à l'ordre 2 de divers $IC$ : $W_s$ (rouge), $RV$ (bleu), Jeffreys (gris), $X + U$ (orange ), Wilson (vert), $AC$ (rose), $BO$ (noir).

L'approximation du biais de différents intervalles de confiance, en fonction de $p$ [Fig. 3] montre que l'intervalle de Wald (en rouge) sous-estime le taux nominal de recouvrement (95%), en particulier lorsque $p$ se rapproche de 0 ou de 1. L'intervalle de Wilson (en vert) ne présente (quasiment) pas de biais, alors que l'intervalle d'Agresti-Coull est légèrement conservatif (la probabilité réelle de recouvrement dépasse légèrement 95%lorsque $p$ n'est pas trop voisin de 0 ou de 1). De la même manière, on pourrait représenter l'amplitude moyenne en $p$ des oscillations en fonction de $n$ et observer une décroissance lente de cette quantité, elle même oscillante pour tous les $IC$ considérés. Cette amplitude est cependant plus élevée pour l'intervalle de Wald standard.

Il convient de confronter également la précision de ces intervalles, de la même manière que l'on compare la puissance de tests d'hypothèse. Bien que le critère de l'espérance de la longueur ait été critiqué (Lehmann, 1959), nous l'utilisons pour sa simplicité.

Il est défini par : $EL_p(IC_*) = E_p[u_*(\hat{p}, n) - l_*(\hat{p}, n)]$

On observe [Fig. 4] que l'intervalle de Wald (en rouge), quoique laxiste, est d'une médiocre précision en regard des intervalles de Wilson(en vert) ou d'Agresti-Coull (en rose). La précision de ce dernier se dégrade cependant nettement lorsque $p$ se rapproche de 0 ou 1. Leur précision est comparée à celle d'autres intervalles que l'on trouve dans la littérature (l'intervalle du rapport de vraisemblance ($RV$), l'intervalle Bayésien avec a priori de Jeffreys), et trois autres intervalles que nous allons également présenter : les intervalles de Clopper-Pearson ($CP$), $Mid$-$P$ et $X+U$.

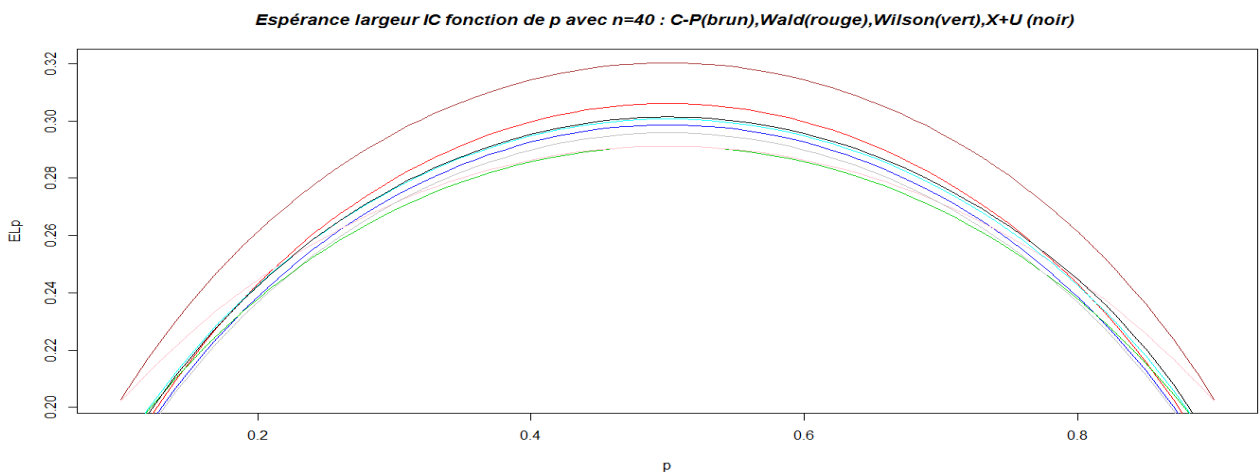

Fig.4 : Espérance de la longueur des divers $IC$ considérés : Wilson (vert), $AC$(rose), $X + U$ (noir),Jeffreys (gris), $RV$ (bleu), Mid-P(cyan), X+U (noir),Wald (rouge), Clopper-Pearson (brun).

## 2. L'intervalle Mid-P et l'intervalle de Clopper-Pearson

Pour éviter d'avoir recours à une approximation, les ouvrages de statistique plus avancés, en particulier ceux orientés vers les applications industrielles (cf. normes AFNOR), recommandent l'intervalle « exact » de Clopper-Pearson $CP$. Cet intervalle est basé, comme l'intervalle de Wilson, sur l'inversion de la famille des régions d'acceptation des tests $H_0 : p = p_0$, $p_0 \in (0,1)$, mais sans avoir recours à l'approximation normale comme ce dernier. Les limites de confiance de $IC_{CP}$ sont définies par les valeurs de $p_0$ solutions de :

$$\sum_{k=x}^{n} \binom{n}{k} p_0^k (1-p_0)^{n-k} = \alpha/2 \quad \text{et} \quad \sum_{k=0}^{x} \binom{n}{k} p_0^k (1-p_0)^{n-k} = \alpha/2 \text{ où } x = \sum_{i=1}^{n} x_i$$

avec pour borne inférieure 0 lorsque $x = 0$ et pour borne supérieure 1 lorsque $x = n$. Cet intervalle a la propriété de garantir une probabilité de recouvrement au moins égale à sa valeur nominale[Fig. 5]. En contrepartie, il est bien trop conservatif et donc peu précis [Fig.4, en brun]. Cet intervalle n'est pas un bon choix en pratique.

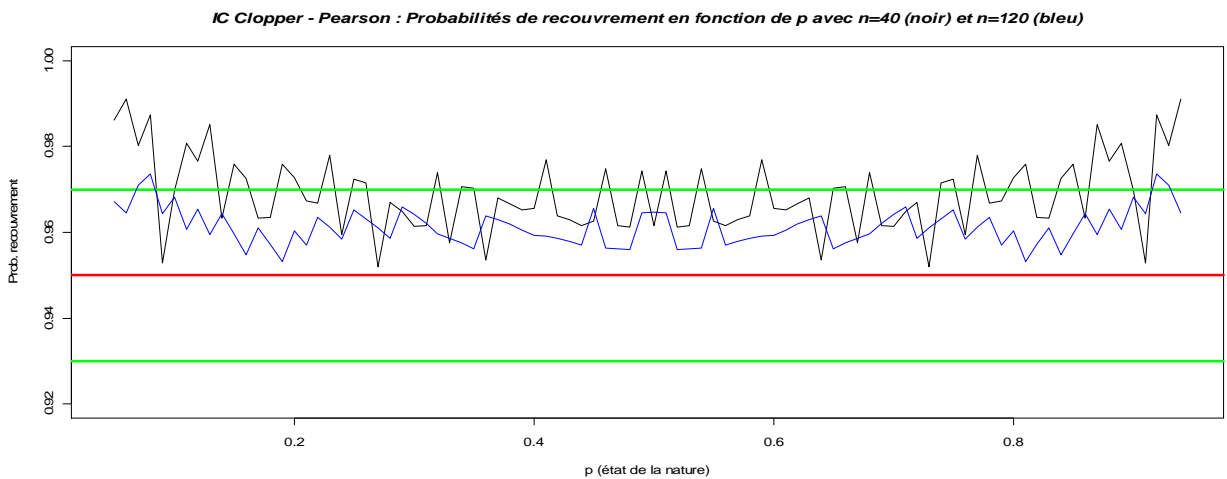

Fig. 5 : Probabilité exacte de recouvrement de $p$ par l'intervalle de Clopper-Pearson, en fonction de $p$ ; $n=40$ (noir), $n=120$ (bleu).

Les limites de l'intervalle de Clopper-Pearson sont obtenues par inversion des inégalités $\Pr_p(X \le l(p)) \le \alpha/2$ et $\Pr_p(X \ge u(p)) \le \alpha/2$ où $u(p)$ et $l(p)$ sont les limites de l'intervalle de fluctuations de niveau $1 - \alpha$. Les limites de confiance $u_{CP}(x)$ et $l_{CP}(x)$ satisfont les inégalités $\Pr_{u_{CP}}(X \le x) \le \alpha/2$ et $\Pr_{l_{CP}}(X \ge x) \le \alpha/2$. Pour éliminer ce conservatisme, on inverse non pas $\Pr_p(X \le l(p)) \le \alpha/2$ mais $\Pr_p(X \le l(p)-1) + (1/2)\Pr_p(X = l(p)) \le \alpha/2$ pour obtenir $u_{MidP}(x)$, pour cela, on résout : $\Pr_{u_{MidP}(x)}(X \le x-1) + (1/2)\Pr_{u_{MidP}(x)}(X = x) = \alpha/2$. De manière symétrique, on obtient $l_{MidP}(x)$. Bien entendu, on perd la garantie du niveau nominal de recouvrement.

## 3. L'intervalle $X + U$ de Stevens.

Cet intervalle est un peu oublié, peut-être à tort. Il est fondé sur la statistique $X + U$ où $X = \sum X_i$ et $U$ la loi uniforme sur $[0,1]$ indépendante de $X$. Autrement dit, aux données sont ajoutées un aléa indépendant des données. La loi de $X + U$ est continue et à rapport de vraisemblance monotone. On en déduit (Lehmann, chap.3 et 5) qu'il existe un intervalle de confiance de niveau exact $1 - \alpha$ pour $p$ appelé parfois intervalle de Stevens. Il a de bonnes propriétés :

(i)  il est sans biais : $\Pr_p \left( p \in IC_{X+U} \right) \geq \Pr_p \left( p' \in IC_{X+U} \right), \forall p' \neq p$

(ii) il est uniformément le plus précis dans la classe des intervalles fondés sur la statistique $X + U$.
L'aléa représenté par $U$ est d'une certaine manière une substitution de l'aléa sur la probabilité exacte des intervalles présentés ci-dessus (E.Pearson). On observe [Fig.4] que sa précision au sens de l'espérance de la longueur de l'intervalle demeure correcte même pour $n$ modéré ($n = 40$), bien meilleure que celle des intervalles de Wald ou de Clopper-Pearson, identique à celle de Mid-P, moins bonne cependant que celle de l'intervalle d'Agresti-Coull. En contrepartie, la probabilité de recouvrement nominale est garantie, ce qui n'est pas le cas des autres intervalles (hormis $C$-$P$).

Nous pensons que l'intervalle de Wald standard ne devrait plus être enseigné, car hormis pour de grands échantillons comme ceux que l'on trouve par exemple dans les sondages à l'échelle nationale, il est trop laxiste. Une simple règle (+2 succès, +2 échecs) donne l'intervalle d'Agresti-Coull, proche de l'intervalle du score de Wilson pris comme référence par la plupart des auteurs. L'expression de ce dernier est malheureusement difficile à interpréter. Nous recommandons l'intervalle d'Agresti-Coull. L'intervalle de Clopper-Pearson, qui n'est pas fondé sur une approximation normale, a pour seul avantage de garantir le niveau de recouvrement nominal. On doit lui préférer l'intervalle de Stevens.